\documentclass[12pt]{article}
\usepackage[utf8x]{inputenc}
\usepackage{amssymb}
\usepackage{color}

\newtheorem{Theorem}{Theorem}
\newtheorem{Lemma}{Lemma}

\begin{document}

\title{
On rate of convergence 
for infinite server Erlang--Sevastyanov's problem 
}


\author{A. Yu. Veretennikov\footnote{School of
Mathematics, University of Leeds, LS2 9JT, Leeds, UK, \& Institute
of Information Transmission Problems, B. Karetny 19, 127994,
Moscow, Russia,
email: a.veretennikov @ leeds.ac.uk}\hspace{1mm}
}


\maketitle

\begin{abstract}
Polynomial convergence rates in total variation 
are established in Erlang--Sevastyanov's type problem with an infinite number
of servers and a general distribution of service under 
assumptions on the intensity of serving. 
\end{abstract}

\section{Introduction}
\label{intro}
The formulae by Erlang provided explicit expressions for percentages 
of lost customers in certain queueing systems in the stationary regime 
\cite{Erlang1917}. 
Erlang models still remain highly important in the modern world. However,  
what is crucial for applications and what is lacking in Erlang's old 
results and some further studies is a knowledge of rate of convergence 
to a stationary regime. In fact, this ``extended'' Erlang's problem -- with 
estimated convergence rates -- is not fully solved even nowadays. For a 
long time, estimations of convergence rate (mostly of exponential 
decay) were known only for the cases where service times 
have exponential distributions and under some additional assumptions,  
cf. \cite{Asmussen}, \cite{Zei}, et al. 
Bounds for the rates of convergence to stationary regimes for 
close systems -- but not precisely Erlang's ones -- were a subject of study in many papers, 
see below. 
It is widely accepted that any important characteristic of quality of
any queueing system in practice is computed in a stationary regime and 
it is, of course, a rare case where this characteristics is available 
in a more or less explicit form, cf., for example, \cite{Pechinkin1979}. 
However, if the rate of convergence 
is unknown, then the error is unknown either. 
Modelling may be some alternative 
to theoretical bounds, nevertheless, it cannot fully replace a 
rigorous theoretical analysis.

Our main goal is to attack the general non-Markov case with 
non-exponential service times for classical telephone systems. 
The key system to be studied is similar to one investigated in 50s by 
Sevastyanov and in the following three decades by other 
researchers. 
This system consists of a finite (as in 
Sevastyanov's works), or infinite 
(as in some other works) number of servers; 
the incoming flow of customers is a conditional Poisson process with 
intensity that may depend on the number of customers in the system; 
in particular, it may increase linearly 
if this number is large. Each customer 
upon his arrival goes to one of the free servers, or -- in the finite 
case -- it may be lost if all servers are busy. All service times are
random variables independent of each other and all have the same
distribution function. Such models even with a finite number of servers 
usually do not satisfy conditions of Doeblin--Doob's ergodic theorem about 
{\em uniform convergence} 
\cite{Doeblin}, \cite[Ch.5-6]{Doob}.

The celebrated B. A. Sevastyanov's
ergodic theorem \cite{Seva} (see also \cite{Seva2}) for Markov processes in general state
spaces provided for the first time not only existence and
uniqueness of stationary distribution for ``telephone systems'', but also
convergence in total variation. This was a pioneering result where such
convergence is {\em non-uniform} with respect to the initial data or
distribution and does not follow directly from Doeblin--Doob's ``uniform''
ergodic theory. The corollary of Sevastyanov's ergodic 
theorem for queueing (``telephone'') model 
will be briefly recalled below. 
Practically simultaneously with \cite{Seva2}--\cite{Seva}, T. E. Harris \cite{Harris1956} suggested his  
method to study stationary measures of recurrent Markov processes; a
presentation of his results and ideas, as well as of their further 
development -- including studies of convergence rates -- 
may be
found in \cite{Baxendale2011}. It may be noted that one of the basic ideas of
this theory -- to exploit moments of ``regeneration'' of the process -- was
proposed in the fourtees in \cite{Doeblin} and further devepoled in
\cite{Kolmogorov} in relation to a very close issue of local limit theorems, 
which may serve as a background for coupling. 
A few years earlier than Sevastyanov and shortly after \cite{Kolmogorov}, 
Fortet proved \cite{Fortet50} that a
stationary distribution exists in ``Sevastyanov's case'' under a bit stronger 
assumptions than eventually in \cite{Seva} (existence of a density was assumed),
along with the form of this distribution; however, he did not study uniqueness 
nor convergence. Some special important part of the main result of \cite{Seva}
-- related to the property of ``insensitivity'' (see below) -- was
also rigorously obtained in \cite{Khinchin62}. The latter paper was 
published only in 1963, 
however, as quite reasonably suggested by the Editor   
(B. V. Gnedenko) 
of the volume of A. Ya. Khinchin's works 
in \cite[The Editor's Introduction]{Khinchin63}, 
the paper was, in fact, fully prepared to publication in 1954--1955.
Earlier, the original Erlang's formulae with
exponential service time distribution were extended
on systems with an infinite number of
servers \cite{Khinchin55}, and later (1965) this result was tackled by a
different method in \cite{Mejzler1965}. Among all these results,  \cite{Seva}
remains the most advanced achievement 
in that period. 

More general systems -- with infinitely many  
servers and/or with more involved disciplines of serving -- were 
studied further in \cite{Falin88}, \cite{Falin89}, \cite{Ivnickii1982},
\cite{Koenig1975},
\cite{Kov1962}, \cite{Matthes1962},  \cite{Schassberger1977},  
\cite{Ve1977}, \cite{Ve1981}, 
et al. Even quite recently, results in 
this direction were still under investigation under the name of 
``insensitivity'' of a stationary regime (i.e., where there are some 
general invariants of a stationary distribution, which depend 
on the service time distribution only through its 
mean value) for advanced versions of Erlang type models in 
\cite{AdanWeiss2011}, 
\cite{Bonald2002}, 
\cite{Massoulie2007}, 
\cite{Walton2011}, 
\cite{Zachary2007}. 
Note that most of these papers -- with the exception of 
\cite{Ivnickii1982} and \cite{Ve1981}  -- 
do not cite two other pioneering 
publications \cite{Kov1962}--\cite{Kov1962b} and none of them except 
\cite{Falin88} tackles convergence rates; in the latter paper, the result 
about convergence rate bounds could be called partial in comparison to our 
Theorem \ref{Thm1} below. 

Sevastyanov's version of ergodic theorem 
\cite{Seva2}--\cite{Seva}  also proved to be
useful in some extensions, in particular, in the case $N=\infty$, see
\cite{Ve1977}. For other versions of such extensions see
\cite{Matthes1962},
\cite{Schassberger1977}, 
et al. (regrettably, the former publication \cite{Matthes1962}  
is still hardly easily available even nowadays).

{\em Exponential} convergence rate for infinite server  
systems of Sevastyanov's type (and a little more general) 
with {\em non--exponential} service time distributions may be found in 
\cite{KV1997}; however, the method used there was not suitable for 
weaker sub--exponential rates under weaker assumptions. 
Establishing such weaker convergence rates for a wider class of 
queueing systems of Sevastyanov's type is the main 
goal of the paper. 

It would be an extremely hard task to mention all important publications
where convergence rates for general Markov processes -- or, indeed, just for
applications to queueing models -- were studied; some of them may be found
among the references below, or in the literature provided in these references. A
very incomplete list of names of major contributors includes Kalashnikov
\cite{Kalashnikov}, \cite{Kalashnikov93}, Borovkov
\cite{Borovkov84}, \cite{Borovkov},  
Tuominen and
Tweedie \cite{TT79}, \cite{Tweedie1983}, Thorisson \cite{Thor85}, \cite{Thor00}, et
al. 
Results about convergence rates 
close to the Theorem \ref{Thm1} below for similar but yet a 
bit different systems may be found in the fundamental monograph \cite{Thor00}, 
where, in particular, the Theorem 7.2 establishes convergence in total
variation, 
\begin{equation}\label{noest}
\varphi(t)\; \|\mu_t - \mu\|_{TV} \to 0, \quad t\to \infty.
\end{equation}
Recall that the total variation distance between two
probability measures on a measurable space $(\Omega, {\mathcal F})$ is
defined as
to 
\[
 \|\mu - \nu\|_{TV} : = 2\,\sup_{A\in {\mathcal F}}\, |\mu(A) - \nu(A)|. 
\]
For certain more particular models see also \cite{Sigman92} and
\cite{TT79}. 
The background idea of  the approach in \cite{Thor00} 
is to use estimates of the 
rate of convergence in the law of large numbers (LLN); its 
implementation is involved and uses regeneration 
technique. In our model with an infinite number of servers ($N=\infty$) 
it is unclear how to use LLN directly and we use another method (eventually 
leading to LLN, too) based on a ``markovisation'' of the
system -- due to Sevastyanov -- and on a local ``infinitesimal'' condition 
on the basis of {\em service intensity}, $h(t)$, which allows to construct
Lyapunov functions. Regeneration is also in use in this paper, which gives 
more precise {\em bounds} for the distance in
(\ref{noest}) and continues studies of various rates of convergence and mixing for a variety of  
Erlang--Sevastyanov's type models commenced recently in 
\cite{V2009} -- \cite{V2013}. The model in this paper is  {\em
non--Markov}.

The paper is arranged as follows. Section one is introduction. 
Section two contains the setting and a brief reminder of 
Sevastyanov's result. Section three is devoted to the main result 
of this paper -- polynomial convergence -- Section 
four to some auxiliaries and Section five to the
proof of the main result.

\section{The setting: Erlang -- Sevastyanov system}
We consider the model with 
$1\le N\le \infty$ (including $\infty$) identical servers working independently, 
with a distribution function $G$ of service time. The incoming flow of
customers is conditionally Poisson with intensity
$\lambda_n$ given that $n$ servers are busy at the
 moment ($0\le n\le N$).  All service times on all servers are independent 
 on each other and on the incoming flow.
A newly arrived customer
chooses any server which is not busy and its serving immediately starts. 
If $N<\infty$ and all servers are busy, then a new customer is lost or
blocked; if $N=\infty$, then under reasonable assumptions the number of
customers is finite at all times and no loss is possible.  
A customer which was served, immediately quits the system. We assume that at any
moment $t$ the elapsed service times of all customers in the system, say,
$X^1_t,
\ldots, X^n_t$ are known; the process $X_t = (X^1_t, \ldots, X^n_t), \, t\ge 0$
is Markov (cf. the Lemma \ref{Le1} below); if there is no customers at $t$, then we denote
$X_t=\Delta_0$ (note that $X_t=0$ and $X_t=\Delta_0$ have different meanings).
At $t=0$, only finitely many servers may be busy. 
Following \cite{Seva}, we assume that 
a newly arrived customer is assigned a coordinate $X^k=0$
with any $k=0, \ldots, n+1$ with equal probabilities $(n+1)^{-1}$ if at his
arrival $n$ servers are busy.

The ``non-Markov property'' of this system signifies that the number $n=n_t \in
Z_+$ of
customers at any time $t$ is, generally speaking, {\em not} a Markov process (of course, unless the {\em intensity} of serving only depends on $n_t$). 
However, we make it Markov by considering it as $(X_t)$ in the
following extended state space ${\mathcal X}$ of a variable
dimension, as in \cite{Seva}: 
${\mathcal X}$ is a union of finitely many (if $N<\infty$), or countably many
(if $N=\infty$) subspaces, 
$$
{\mathcal X}_0 = \Delta_0; \qquad
{\mathcal X}_1 = R_+, \quad \ldots, 
\qquad
{\mathcal X}_n = R^n_+, \quad \ldots, 0\le n\le N.
$$
To any $x\in {\mathcal X}_n$ with $n>0$ there correspond $n$ non-negative 
coordinates $(x^1, \ldots, x^n)$, which signify the elapsed time 
of service of any of existing $n$ customers in the system.
%
If there is at least one customer in the system and $x=(x^1,
\ldots, x^n)$ is a vector of the elapsed service times,
then by $n(x)$ we denote this number $n$; if $x = \Delta_0$,
then $n(x):=0$. 

~

\noindent
In \cite{Seva} (see also \cite{Seva2}) it is proved that for $N<\infty$ under
the only condition 
$$
q^{-1} := \int_0^\infty x \,dG(x) < \infty,
$$ 
there is a (unique) stationary distribution $\mu$ with a 
density ($1\le k\le N$), 
$$
p_k(x^1, \ldots, x^k) = p_0\frac{\prod_{i=0}^{k-1}\lambda_i}
{k!} \, \prod_{j=1}^{k} (1-G(x^j)), \quad   
p_0^{-1} = \sum_{j=0}^{N-1}\frac{\prod_{i=0}^{j-1}\lambda_i}
{q^j\, j!},
$$
and, moreover, for any initial distribution the following convergence holds true,
\begin{equation}\label{to}
\|\mu_t - \mu\|_{TV} \to 0, \qquad t\to\infty,
\end{equation}
where $\mu_t$ is a marginal distribution of the (Markov) process at $t$; 
below $\mu_t^x$ will stand for the marginal distribution given initial state $x$. 

~  

\noindent 
{\it Remark 1.} 
In \cite{Ve1977} this result was extended to  $N=\infty$ 
under the condition $q^{-1} \limsup_{n\to\infty}\lambda_n/(n+1)<1$. 
By integrating $\int p_k(x)\,dx$, we obtain the
stationary 
probabilities $p_k$ of $k$ customers in the system, which depend on $G$ 
only through its mean value; this property is called ``insensitivity'',
$$
p_k = p_0\frac{\prod_{i=0}^{k-1}\lambda_i}
{q^{k} k!}, \quad 1\le k\le N,
$$
and it is an object of studies for various queueing models  until nowadays.

\section{Convergence rate bounds: Main result ($N=\infty$) }
We consider Erlang--Sevastyanov's system with $N=\infty$. 
The {\em service intensity} $h(t)$ is defined as follows, which 
we will assume to be {\em bounded} ($h\in B(R_+)$) 
$$
h(t) : = \frac{g(t)}{1-G(t)}, \quad t \ge 0, \quad \mbox{where}\;\;
g(t)  = G'(t).
$$
Notice that $h \equiv \mathop{\mbox{const}}$ means that the service time 
has an exponential distribution, in which case 
(and in a bit more general one) a sufficient condition for exponentially fast 
convergence to the stationary distribution has been established in
\cite{KV1997}. 
%
In all cases,  $\displaystyle 1-G(t)=\exp\left(-\int_0^t h(s)\,ds\right)$.
Let for $a,m>0$, 
$$
L_{m,a}(x):= 
\left(\sum_{j=1}^{n(x)}(1+x^j)^{m}\right)^{a} \;\; (x\not= \Delta_0), 
\quad \& \quad L_{m,a}(\Delta_0):= 0. 
$$
To avoid triviality due to a degeneracy, we impose a condition
\begin{equation}\label{nondeg}
 \lambda_0>0.
\end{equation}

\begin{Theorem}\label{Thm1}
Assume (\ref{nondeg}), 
$h\in B(R_+)$ 
and  existence of $C_0 >0$ 
such that 
\begin{equation}\label{ash}
h(t) \ge  \frac{C_0}{1+t}, \quad t\ge 0,
\end{equation}
and 
\begin{equation}\label{aslambda}
C_0 - 2(1 + 2\Lambda) > 0,   
\end{equation}
where
\begin{equation}\label{Lambda}
\Lambda : = \sup_{n\ge 1} \left(\frac{\lambda_n}{n}\right) < \infty.
\end{equation}
Then for any $k>0$ small enough
there exist real values $C>0$, $a>1$ and $m > 1$ such that 
for any $X_0=x$, 
\begin{equation}\label{main}
 \|\mu^x_t - \mu \|_{TV} \le 
 C \,\frac{1+L_{m,a}(x)}{(1+t)^{k+1}}.
\end{equation}
For $\Lambda$ fixed, $k$ may be chosen large if $C_0$ is large enough, see
(\ref{am}) below. 
\end{Theorem}

\noindent
{\em Remark 2.} Without (\ref{nondeg}) -- i.e., for $\lambda_0=0$ -- the Theorem \ref{Thm1} is still valid
with a trivial stationary distribution, $\delta_{\Delta_0}$ and it does follow
from the proof below with minimal changes. 
 
~

\noindent
{\em Remark 3.} 
As we shall see in the proof, for a substantial part of the proof it suffices to assume a slightly weaker assumption
\begin{equation}\label{aslambda0}
C_0 - (1 + \Lambda) > 0.
\end{equation}
However, in the end of the calculus the full 
version of (\ref{aslambda}) will be needed. More precisely, 
we will, actually, use
\begin{equation}\label{am}
C_0 > 
(a+(k\vee 1)/m)(m + \Lambda
\,2^{a-1+(k\vee 1)/m}). 
\end{equation}
The latter bound is available with {\em some} $a,m>1$, at least, for small $k>0$ -- in fact, for $k\le m$ -- 
under the assumption (\ref{aslambda}). 
As one more example, with $m=k$, the latter sufficient condition only in terms
of $C_0, \Lambda$ and $k$ reads, 
\[
 C_0> 2(k+2\Lambda),
\]
as in this case there exists $a>1$ for which (\ref{am}) holds.  
Also notice that large values of $k$ in (\ref{am}) are available under
$C_0$ large
enough, or under just $C_0>1$ but with $\Lambda$ small enough, which agrees with the intuitive
idea that stability is stronger if intensity of service is in some
sense
significantly greater than intensity of arrivals. However, emphasize that $C_0$
itself is {\em not} intensity of serving itself, but only a multiplier in a lower bound 
for this function (\ref{ash}). 

~

\noindent
{\em Remark 4.} Of course, the greater $C_0$, the more moments has the distribution of serving $G$. However, the method requires existence of intensity $h$. It would be interesting to relax the latter assumption.

\section{Construction, martingales, estimates,  strong Markov property}
We will use notations $x=(x_1, \ldots, x_n)$ and $X=(n,x)$ -- and also $E_X\equiv E_x$ -- and for any such 
$x\in {\mathcal X}_n$ with $n\ge 1$ define with any $1\le j\le n$,
\[
 x^{(+j)}:=(x_1, \ldots, x_j, 0, x_{j+1}, \ldots, x_n), \;\; 
\& \;\; 
 x^{(-j)} := (x_1, \ldots, x_{j-1}, x_{j+1}, \ldots, x_n).
\]
To work with Lyapunov functions, it is very useful -- if not compulsory -- to know 
that the process is {\em strong Markov}. In continuous time it is not 
automatic and should be justified. About {\em Markov} property for finite $N$, cf. 
\cite{Fortet50} and \cite{Seva}.

\noindent
{\em Some preliminaries: generators and martingales.} Suppose for a while that  $h\in C_b(R_+)$; a bit later we will show how to relax this assymption. 
Before the next Lemmae  
we recall some well-known links between Markov processes and martingales, 
which seem to be a bit less popular language in queueing theory (e.g., cf. \cite{Borovkov84}). 
The generator (infinitesimal operator) ${\mathcal G}$ of the process $X$ in 
the space ${\mathcal X}$ with a Borel topology ${\mathcal B}$ on all 
subspaces ${\mathcal X}_n$ (with a convention that ${\mathcal X}_n$ is open and closed 
for each $n$) and $\sup$-norm for $C({\mathcal X}, {\mathcal B})$
is an operator ${\mathcal G}$ such that (see \cite{Dynkin})
\begin{equation}\label{gen}
\sup_{X\in {\mathcal X}}\left|\frac{E_X f(X_t)-f(X)}{t} 
- {\mathcal G}f(X)\right| \to 0, \quad t\to 0,
\end{equation}
for all $f$ from the {\em domain} ${\mathcal D}_{\mathcal G}$ of ${\mathcal G}$, 
which is usually a hard task to determine precisely and which is usually enough to 
have a wide enough sub-class of. 
In our case, 
it follows from (9)--(10) and continuity of $h$ that for 
$f\in C^1_0({\mathcal X})$ -- with one continuous derivative and compact support 
-- i.e., $f(X)$ vanishes if $n\ge N$ or if $\sup_i x^i\ge N$  for some $N$ 
-- (\ref{gen}) holds with 
\begin{eqnarray}\label{ext-gen}
{\mathcal G}f(X) = 
\sum_{i=1}^{n(X)}\left(\frac{\lambda_{n(X)}}{n(X)}(f(X^{(+i)}) - f(X)) 
+  h(X^i) (f(X^{(-i)}) - f(X)) \right.
 \nonumber \\ \\ \nonumber
\left. 
+ \frac{\partial}{\partial x^i}f(X)\right) 
\times 1(n(X)>0) + \lambda_0 (f(0) - f(X))\, 1(n(X)=0). \hspace{1cm}
\end{eqnarray}
By Dynkin's formula \cite[see, e.g., corollary from the formula (1.36) as $\lambda\to\infty$]{Dynkin}, 
\begin{equation}\label{dynkin1}
 E_{X_0} f(X_{t}) - f(X_0) = E \int_0^{t} 
 {\mathcal G}f(X_s)\,ds
\end{equation}
for any $f\in C^1_0$. For functions $(f(t,X), \, t\ge 0, \, X\in {\mathcal X})$ 
of class $C^1_0$ with respect to $(t,X)$ -- which vanish for large $n(x)$ 
and for large $X\in {\mathcal X}_n$ for any fixed $n$ -- Dynkin's formula for the process 
$(t,X_t)$ reads,  
\begin{equation}\label{dynkin2}
 E_{X_0} f(t, X_{t}) - f(0,X_0) = E \int_0^{t} 
 \left(\frac{\partial}{\partial s}\, f(s,X_s) + {\mathcal
G}f(s,X_s)\right)\,ds. 
\end{equation}
Note that, at least, intuitively this equality as well as  (\ref{dynkin1}) may be regarded as a  
complete probability formula, 
as the right hand side presents all possible developments of the trajectory from $0$ to $t$.  
In terms of martingales (cf. \cite{LSh}), (\ref{dynkin2})  is equivalent to 
saying that the process
\begin{equation}\label{mart1}
M_{t}=f(t, X_{t}) - f(0,X_0) - \int_0^{t} 
 \left(\frac{\partial}{\partial s}\, f(s,X_s) + {\mathcal
G}f(s,X_s)\right)\,ds,\; \;t\ge 0
\end{equation}
is a martingale.  
Recall for the sequel that a process is called a {\em local} martingale if there exists a sequence of
stopping times $\tau_n\to\infty$ a.s. such that the stopped process 
\(
 M_{t\wedge\tau_n}
\)
is a martingale for each $n$ (cf. \cite{LSh}). In turn, the statement
(\ref{mart1}) may be
equivalently (by definition) rewritten in the differential form as 
\begin{equation}\label{difform}
 df(t,X_t) = \left(\frac{\partial}{\partial t}\, f(t,X_t) + {\mathcal
G}f(t,X_t)\right)\,dt + dM_t
\end{equation}
for any $f\in C^1_0$ (${\mathcal G}$ is defined above). The latter formula itself 
is {\em local} -- i.e., written at 
any given $(t, X_t)$ and its small neighbourhood -- so it may be extended from $C^1_0$ to $C^1$, 
under a natural convention that the process is stopped on the exit 
from some neighbourhood of the state $X_t$; of course, this may require 
a localizing stopping time procedure if using the integral form (\ref{mart1}) and possibly some justification that $M$ is a martingale (and not just a local martingale). 

All the above starting from the formula (\ref{gen}) 
work well if the intensities are continuous. If this is wrong, the limit in (\ref{gen}) may just not exist. Nevertheless, following \cite{Davis} or \cite{GK}
it is possible to define the process by using a notion of {\em extended generator}, that is, an operator for which Dynkin's formulae (\ref{dynkin1})  and (\ref{dynkin2}) hold. The action of extended generator on functions is given by the same expression in (\ref{ext-gen}).

\begin{Lemma}\label{Le1}
Under the assumptions (\ref{Lambda}) and $h\in B(R_+)$  
the process $X$ exists, has a unique distribution and is Markov and strong Markov. The Dynkin's formulae (\ref{dynkin1})  and (\ref{dynkin2}) hold for any 
$f(x)\in C^1_0$ and $f(t,x)\in C^1_0$. 
\end{Lemma}
{\bf Proof.} Existence (for possibly discontinuous $h$) follows from the results on piecewise linear or piecewise deterministic Markov processes in \cite{GK}, \cite{Kalashnikov77}, \cite{Davis}, as well as do uniqueness and Markov and strong Markov properties. The non-explosion is implied by the condition  (\ref{Lambda}), for example, due to 
\cite[Ch. 1.3.3]{GK}. Both Dynkin's formulae follow from \cite{Davis}. Another way to show Dynkin's formula for a sightly different model was suggested in \cite{VZ1}.  
The Lemma \ref{Le1} is proved. 

~

We admit the following convention for stochastic
differentials:
\[
 A_t\,dt + dM_t \le B_t\,dt + dM_t \qquad \mbox{iff} \qquad  
 A_t \le B_t, \quad \forall \; t\ge 0, \;\;\mbox{a.s.}
\]
Recall that the process is called {\em cadlag} iff it is right 
continuous with left limits at any $t$. 
Note that $M_t$ in (\ref{mart1}) is cadlag because the right hand side is.
Below we use convention 
$n^{-1}\sum_{i=1}^na_i \equiv 1$ for any real values $(a_i)$ if $n=0$. 

\begin{Lemma}\label{Le3}
Under the assumptions (\ref{Lambda}) and $h\in C_b(R_+)$, 
\begin{eqnarray}\label{M2}
L_{m,a}(X_t) - L_{m,a}(X_0) 
= \int_0^t \lambda_{n(X_s)} \left(\frac1{n(X_s)} \sum_{i=1}^{n(X_s)} (L_{m,a}(X^{(+i)}_s) -
L_{m,a}(X_s))\right) \,ds
  \nonumber \\ \\ \nonumber 
+ \int_0^t \left[\sum_{i=1}^{n(X_s)} h(X^i_s) \left(L_{m,a}(X^{(-i)}_s) - L_{m,a}(X_s)\right) 
+ \sum_{i=1}^{n(X_s)} \frac{\partial}{\partial x^i}L_{m,a}(X_s)\right]\,ds +M_t,
 \end{eqnarray}
with some martingale $M_t$. 
\end{Lemma}
{\bf Proof}
follows from (\ref{difform}) with $f(t,X)\equiv L_{m,a}(X)$ (see the remark above 
about extension of (\ref{difform}) to $C^1$), due to 
\begin{eqnarray}\label{MM}
dL_{m,a}(X_t) 
 = {\mathcal G}(X_{t})\,dt + dM_t 
= \lambda_n \left(\frac1n \sum_{i=1}^{n} L_{m,a}(X^{(+i)}_t) -
L_{m,a}(X_t)\right) \,dt
  \nonumber \\  \\ \nonumber 
+ \sum_{i=1}^n h(X^i_t) \left(L_{m,a}(X^{(-i)}_t) - L_{m,a}(X_t)\right) \,dt 
+ \sum_{i=1}^{n} \frac{\partial}{\partial x^i}L_{m,a}(X_t)\,dt +dM_t 
  \nonumber \\ \nonumber \\ \nonumber 
\equiv (I_1 - I_2 + I_3)dt+ dM_t,
\hspace{3cm}
 \end{eqnarray}
with $n=n(X_t)$. We shall see below in the Lemma  \ref{Le4} (without a vicious circle) that no localization is needed 
here as all terms in (\ref{M2}) will turn out to be integrable and $M$ is, indeed, a martingale.
The Lemma \ref{Le3} is proved. 

%

\begin{Lemma}\label{Le3a}
Under the assumptions (\ref{Lambda}), $h\in B(R_+)$
and $m,a>1$,  
the following bounds or equalities hold true:
\begin{equation}
I_1 \le 
\Lambda \, a \,2^{a-1} \, L_{m-1,1}(X_{t}) L_{m,a-1}(X_{t}) 
1(X_t \not = \Delta_0)  + \lambda_0 1(X_t = \Delta_0);
\end{equation}
\begin{eqnarray}\label{query}
I_2  
\le \,a\, \|h\|_{B}\,L_{m,a+1}(X_t) 1(X_t\not = \Delta_0);
\end{eqnarray}
\begin{eqnarray}\label{i33}
 I_3 =  a m \, L_{m-1,1}(X_{t}) L_{m,a-1}(X_{t}) 1(X_t \not= \Delta_0); 
\end{eqnarray}
\begin{eqnarray}\label{momm}
 E_x L_{m,a} (X_{t}) 
\le (L_{m,a} (x)+\lambda_0 t)  \exp((\Lambda a 2^{a-1} + m a)t). 
\end{eqnarray}
If in addition 
\begin{equation}\label{asam0}
 C_0 > a(m + \Lambda \,2^{a-1}), 
\end{equation}
then also 
\begin{eqnarray}\label{query2}
I_2\ge 1(X_t\not = \Delta_0) \,C_0 L_{m-1,1}(X_{t})
L_{m,a-1}(X_{t}).
\end{eqnarray}

\end{Lemma}

\noindent
{\bf Proof.}
Let us establish the bound for $I_1$. 
Notice that for $y =  x +1 \ge 2$ and $a>1$ 
we have $y^{a-1}= (x+1)^{a-1}\le (2x)^{a-1}$, and, hence,  
\(
y^a-x^a \le a(y-x)y^{a-1} \le a\,2^{a-1}\,x^{a-1}.
\)
Indeed, the first bound here follows for $y\ge x>0$ and $a>0$ from  
\[
 \frac{d}{dx} (y^a-x^a) = -ax^{a-1} \ge -ay^{a-1} =  
 \frac{d}{dx} (a(y-x)y^{a-1}).
\]
So, we estimate, 
\begin{eqnarray*}
\frac1n \sum_{i=1}^{n} L_{m,a}(X^{(+i)}_t) -
L_{m,a}(X_t) \hspace{2cm}
 \\\\
= \left( 
\left((1+0)^m +\sum_{j=1}^{n}(1+X_t^j)^{m}\right)^{a} 
- \left(\sum_{j=1}^{n}(1+X_t^j)^{m}\right)^{a}
\right)
 \\\\
\le a\,2^{a-1}\, \left(\sum_{j=1}^{n}(1+X_t^j)^{m}\right)^{a-1}
 = a\,2^{a-1}\,L_{m,a-1}(X_t).
\end{eqnarray*}
Hence, due to the inequlaity $n(X_t)  \le L_{m-1,1}(X_{t})$, we get,  
\begin{eqnarray}\label{i1}
 I_1 
= \lambda_n \left(L_{m,a}(X'_t) - L_{m,a}(X_t)\right) 
\le \lambda_n \, a \,2^{a-1} \, L_{m,a-1}(X_{t}) dt
  \nonumber \\ \\ \nonumber
 \le \Lambda n \, a \,2^{a-1} \, L_{m,a-1}(X_{t}) 
 \le \Lambda \, a \,2^{a-1} \, L_{m-1,1}(X_{t}) L_{m,a-1}(X_{t}).
\end{eqnarray}
Further, by taking derivatives, we find,
\begin{eqnarray}\label{i3}
 I_3 = \sum_{i=1}^{n} \frac{\partial}{\partial x^i}L_{m,a}(X_t)
= a \left(\sum_{\ell=1}^{n}(1+X_t^\ell)^{m}\right)^{a-1}
\times \sum_{j=1}^n m \, (1+X_t^j)^{m-1} 
 \nonumber \\ \\ \nonumber 
=  a m \, L_{m-1,1}(X_{t}) L_{m,a-1}(X_{t}). \hspace{3cm}
\end{eqnarray}
The lower bound for $I_2$ under the additional (\ref{asam0}),   
\begin{eqnarray}\label{query1}
I_2 \ge C_0 \sum_{i=1}^n 
(1+X^i_t)^{m-1} 
L_{m,a-1}(X_{t})
= C_0 L_{m-1,1}(X_{t})
L_{m,a-1}(X_{t}).
\end{eqnarray}
Emphasize that $L_{m,a-1}(X_{t})$ stands here in the middle term 
and not $L_{m,a-1}(X^{(-i)}_{t})$ -- the latter would be a little bit
insufficient for our aims -- which is justified in the next few lines.
We used here the elementary inequality for real values $0<x\le y$ and $a>1$, 
\begin{equation}\label{elem}
 y^a-x^a \ge (y-x)y^{a-1}
\end{equation}
(instead of also correct $y^a-x^a \ge a(y-x)x^{a-1}$), 
for $y=\sum_{j=1}^{n}(1+X_t^j)^{m}$ and $x =
\sum_{1\le j\le n, \, j\not=i}^{}(1+X_t^j)^{m}$.
The bound (\ref{elem}) follows from the observation that both sides in
(\ref{elem}) vanish at $y=x$ and the derivative function of the 
right hand side is less than that of the left hand side for $y>x\, (>0)$: 
$$
\frac{d}{dy} \,(y-x)y^{a-1} = ay^{a-1} - (a-1)xy^{a-2} < ay^{a-1} 
= \frac{d}{dy} \,(y^a-x^a).
$$
Hence,  
\begin{eqnarray*}
I_2 = \sum_{i=1}^n h(X^i_t) 
\left( 
\left(\sum_{j=1}^{n}(1+X_t^j)^{m}\right)^{a} 
- \left(\sum_{1\le j\le n, \, j\not=i}^{}(1+X_t^j)^{m}\right)^{a}
\right)
 \\\\
\ge \sum_{i=1}^n 
\frac{C_0}{(1+X^i_t)} 
\left(\sum_{j=1}^{n}(1+X_t^j)^{m} - 
\sum_{1\le j\le n, \, j\not=i}^{}(1+X_t^j)^{m}\right) 
\left(\sum_{j=1}^{n}(1+X_t^j)^{m}\right)^{a-1}
 \\\\
= C_0\, \sum_{i=1}^n (1+X^i_t)^{m-1}  
\left(\sum_{j=1}^{n}(1+X_t^j)^{m}\right)^{a-1}
= C_0 L_{m-1,1}(X_{t})
L_{m,a-1}(X_{t}).
\end{eqnarray*}
The upper bound for $I_2$ follows from its definition and from the 
remark that $n(x) \le L_{m,1}(x)$ and $L_{m,1}L_{m,a}=L_{m,a+1}$. 

~

\noindent
Further, for any $t\ge 0$,
\begin{eqnarray*}
I_1 -I_2+ I_3 \le  
(\Lambda a 2^{a-1} + m a) \, L_{m-1,1}(X_{t}) L_{m,a-1}(X_{t}) + \lambda_0.
\end{eqnarray*}
So, from (\ref{MM}) and by virtue of Fatou's lemma -- and using a localization 
for $M$ if needed so as to
vanish expectation of the (local) martingale term -- we get, 
\begin{eqnarray*}
E_x L_{m,a} (X_{t\wedge \tau_n}) 
\le L_{m,a} (x) + \lambda_0 t \hspace{2cm}
 \\\\
+ (\Lambda a 2^{a-1} + m a) E_x \int_0^{t \wedge \tau_n} 
L_{m-1,1}(X_{s}) L_{m,a-1}(X_{s}) \, ds
 \\\\
\le L_{m,a} (x)  + \lambda_0 t + (\Lambda a 2^{a-1} + m a) 
E_x \int_0^{t} L_{m,a}(X_{s\wedge \tau_n}) \, ds.
\end{eqnarray*}
By Gronwall's inequality (note that $E_x L_{m,a} (X_{t \wedge \tau_n})$ is bounded), 
\begin{eqnarray*}
E_x L_{m,a} (X_{t \wedge \tau_n}) 
\le (L_{m,a} (x)+\lambda_0 t) \exp((\Lambda a 2^{a-1} + m a)t). 
\end{eqnarray*}
and, as $\tau_n\to \infty$, by Fatou's Lemma, also
\begin{eqnarray*}
 E_x L_{m,a} (X_{t}) 
\le (L_{m,a} (x)+\lambda_0 t)  \exp((\Lambda a 2^{a-1} + m a)t). 
\end{eqnarray*}
The Lemma \ref{Le3a} is proved.

\begin{Lemma}\label{Le4}
Under the assumptions (\ref{Lambda}), $h\in C_b(R_+)$, 
$m,a>1$, 
for any $t>0$, 
\begin{equation}\label{supl}
E_x\sup_{0\le s\le t} L_{m,a}(X_s) 
 <\infty, 
\end{equation}
and the local martingale $M$ in (\ref{M2}) is, in fact, a martingale.
\end{Lemma}
{\bf Proof.}
We estimate, {\em for any $a,m>0$},
\[
E_x\sup_{0\le s\le t} L_{m,a}(X_s) \le L_{m,a}(x) 
+ \int_0^{t} E_x L_{m,a}(X_s)ds + E_x\sup_{0\le s\le t} |M_s|. 
\]
In turn, $E_x\sup_{0\le s\le t} |M_s| \le C_p (E_x |M_{t}|^p)^{1/p}$ 
for any $p>1$ by Doob's inequality (recall that $M$ is cadlag) and further, 

\begin{eqnarray}\label{modm}
|M_t| \le L_{m,a}(X_t) + L_{m,a}(x) 
 + \left|\int_0^t {\mathcal G}(X_{s})\,ds\right| 
\hspace{1.4cm}
 \nonumber \\ \\ \nonumber
\le L_{m,a}(X_t) + L_{m,a}(x)  
+\left|\int_0^t I_1\,ds\right| 
+\left|\int_0^t I_2\,ds\right|
+\left|\int_0^t I_3\,ds\right|.
\end{eqnarray}
So, due to (\ref{momm}), which is valid for any $m, \,a>0$,
\begin{eqnarray*}
E_x |M_{t}|^p \le C_p L_{m,a}(x)^p + C_p E_x L_{m,a}(X_{t})^p +C_{p}  
\int_0^{t} E_x L_{m,a+1}(X_s)^p ds 
 \\\\
=  C_p L_{m,ap}(x) + C_p E_x L_{m,ap}(X_{t}) +C_{p}  \int_0^{t} E L_{m,ap+p}(X_s)\,ds \hspace{0.5cm}
 \\\\
\le C'(L_{m,ap+p}(x)+\lambda_0 t)\exp(Ct)<\infty. 
\hspace{1.5cm}
\end{eqnarray*}
By virtue of H\"older's inequality, this implies (\ref{supl}). The Lemma \ref{Le4} is proved.

\section{Proof of Theorem \ref{Thm1}}
\noindent
{\bf 1.}
Consider a Lyapunov function $L_{m,a}$ at 
any $X_t\not=\Delta_0$ and $m, a > 1$, satisfying also 
(\ref{asam0})
(compare with (\ref{aslambda}) and (\ref{am})). 
The idea of Lyapunov functions in a stochastic 
context is to verify that the `main' negative term prevails 
and that `on average' the process $L_{m,a}(X_t)$ decreases, as long as 
$X_t \not=\Delta_0$. 
From the bounds on $I_1$, $I_2$ and $I_3$ of the Lemma \ref{Le3a} it follows 
that $I_1$ and $I_3$  are {\em dominated} by 
$I_2$. Then it would imply that the stationary measure 
integrates some polynomial. In turn, this would allow to 
extend our Lyapunov function so as to include some 
multiplier that depends on time. The latter would help obtain  
the crucial bound
$
E_x~\tau_0^{k+1}~<~\infty
$
for some $k>0$ (see (\ref{simply}) below)). 
Finally, the latter bound would imply ``coupling'' 
between the original process and its stationary version 
with a certain rate of convergence. 
\noindent
For any $t<\tau_0$ we have,
\begin{eqnarray*}
I_1 - I_2 + I_3 \le - 
(C_0 - \Lambda a 2^{a-1} - m a) \, L_{m-1,1}(X_{t}) L_{m,a-1}(X_{t}) < 0.
\end{eqnarray*}
Denote 
$$
C_{m,\Lambda,a}:= C_0 - \Lambda a 2^{a-1} - m a >0.
$$ 
By Fatou's lemma 
we get, 
\begin{eqnarray}\label{ineqlt}
E_x L_{m,a} (X_{t \wedge \tau_0})  \hspace{3cm} \nonumber \\\\ \nonumber
+ C_{m,\Lambda,a} E_x \int_0^{t \wedge \tau_0} 
L_{m-1,1}(X_{s}) L_{m,a-1}(X_{s}) \, ds
\le L_{m,a} (x),  
\end{eqnarray}
and, as $t\to \infty$,
\begin{eqnarray*}
 E_x L_{m,a} (X_{\tau_0}) 
+ C_{m,\Lambda,a} E_x \int_0^{\tau_0} 
L_{m-1,1}(X_{s}) L_{m,a-1}(X_{s}) \, ds 
\le L_{m,a} (x). 
\hspace{4cm}
\end{eqnarray*}
In particular, $E_x\tau_0 < \infty$ for any $x$ and also 
$E_0\hat\tau_0 < \infty$ with $\hat\tau_0 := 
\inf(t > 0: \, X_t=\Delta_0; \; \exists\, s\in (0,t): \, 
X_s \not=\Delta_0))$. In other words, the process $X$ is positive
recurrent. 
According to the Harris--Khasminskii principle about invariant measures 
(cf., for example, \cite{Ve00}), 
there exists a (unique in our model) invariant measure $\mu$, 
$\mu(A) = c E_0 \int_0^{\hat\tau_0} 1(X_s\in A)\,ds$, which integrates 
the function $L_{m-1,1}(x) L_{m,a-1}(x)$. As noticed by the Referee, under the accepted assumtions both existence and 
uniqueness of this measure also follow straightforward from  \cite{Ve1977}.

~

\noindent
In a moment, we will show
one more elementary inequality 
\begin{equation}\label{gm}
 L_{m,1}(x)^{(m-1)/m} \le L_{m-1, 1}(x), 
\end{equation} 
so that (notice that  $L_{m,a} (x) L_{m,b} (x) = L_{m,a+b} (x)$ and $L_{m,1}(x)^{a} = L_{m,a}(x)$)
$$
L_{m-1,1}(x) L_{m,a-1}(x) \ge L_{m,a-1+\frac{m-1}{m}}(x) = L_{m,a - 1/m}(x),
$$
and
\begin{eqnarray}\label{HH}
 E_x L_{m,a} (X_{t\wedge \tau_0}) 
+ C_{m,\Lambda,a} E_x \int\limits_0^{t\wedge \tau_0} 
L_{m,a-1/m}(X_{s}) \, ds 
\le L_{m,a} (x).
\end{eqnarray}
So, due to Fatou's lemma,
\begin{eqnarray}\label{HH2}
 E_x L_{m,a} (X_{\tau_0}) 
+ C_{m,\Lambda,a} E_x \int_0^{\tau_0} 
L_{m,a-1/m}(X_{s}) \, ds 
\le L_{m,a} (x). 
\end{eqnarray}
Emphasize that both inequalities  
(\ref{HH}) and (\ref{HH2}) have been established under the assumption 
$
C_{m,\Lambda,a} >0,
$ 
that is, 
\begin{equation}\label{encore1}
C_0 - \Lambda a 2^{a-1} - m a >0.
\end{equation}

~

\noindent
{\bf 2.} 
The inequality (\ref{gm}) 
follows 
from the inequalities with $a,b>0, \, \alpha\in (0,1)$ and $c = b/a$ 
$$
(a+b)^\alpha \le a^\alpha +b^\alpha \quad \sim \quad 
(1+c)^\alpha \le 1 +c^\alpha, 
$$
where the latter, in turn, follows from the 
valid inequality 
for the derivatives, 
$$
\alpha(1+c)^{\alpha-1} \le \alpha c^{\alpha-1}. 
$$  

\noindent
{\bf 3.} 
We are now prepared for considering a Lyapunov function 
which depends on time. Let $k>0$ ({\em not} necessarily $k>1$), 
$a, m>0$ and $L_{m,a,k}(t,x): = (1+t)^k L_{m,a}(x)$. 
Similarly to the above and choosing $a, m>1$, we have, 
\begin{eqnarray*}
dL_{m,a,k}(t,X_t)= L_{m,a,k}(t+dt,X_{t+dt}) - L_{m,a,k}(t,X_{t}) 
 \\\\
=  (1+t)^{k}\left[I_1 - I_2 +I_3\right]dt  
+ k(1+t)^{k-1} \,  L_{m,a}(X_{t})\, dt
+ d\tilde M_t
 \\\\
\le -(1+t)^k (C_0 - \Lambda a 2^{a-1} - m a) 
L_{m,a-\frac1{m}}(X_{t}) \, dt  \hspace{0.5cm}
 \\ \\
+ k(1+t)^{k-1} \, L_{m,a}(X_{t})\, dt + d\tilde M_t, \hspace{1.5cm}
\end{eqnarray*}
with some new local martingale $\tilde M$.
Now the task is again to ensure that the negative part in the 
right hand side of the last expression prevails. 
We will be using the 
inequality established in the step 1 above. The second term may 
be split into two parts, 
\begin{eqnarray}\label{twoterms}
I:= k(1+t)^{k-1} \,  L_{m,a}(X_{t})  \hspace{2cm}
 \\\nonumber 
= I\times 1(k(1+t)^{k-1} \, L_{m,a}(X_{t})
\le \epsilon (1+t)^k L_{m,a-1/m}(X_{t}))
  \\\nonumber
+ I\times  1(k(1+t)^{k-1} \, L_{m,a}(X_{t})
> \epsilon (1+t)^k L_{m,a-1/m}(X_{t})).
\end{eqnarray}
The first term here with '$\le \epsilon$', clearly, is dominated by 
the main negative expression if $\epsilon>0$ is small enough,   
$
\epsilon < C_0 - a(m + \Lambda \,2^{a-1}).
$ 
The set of such values $\epsilon$ is non-empty 
as long as $a$ and $m$ are chosen so as to satisfy 
(\ref{asam0}).

  
Let us estimate the second term in (\ref{twoterms}). 
We have, for any $\ell>0$ (later we will choose $\ell = k+\delta$ with small
$\delta>0$),
\begin{eqnarray*}
I\times  1(k(1+t)^{k-1} \, L_{ m,a}(X_{t})
> \epsilon (1+t)^k L_{m,a-1/m}(X_{t}))
 \\\\
\le I \times  \frac{(k \, L_{m,a}(X_{t}))^\ell}
{(\epsilon (1+t) L_{m,a-1/m}(X_{t}))^\ell} 
= I \times  \frac{k^\ell}{(\epsilon (1+t))^\ell} 
L_{m,\ell/m}(X_{t}). 
\end{eqnarray*}
Therefore, the second term  in (\ref{twoterms}) does not exceed
\begin{eqnarray*}
k(1+t)^{k-1} \,  
\times \frac{k^\ell}{(\epsilon (1+t))^\ell}
L_{m,a+\ell/m}(X_t)^{}.
\end{eqnarray*}
Now let us impose conditions on $\ell$: 
let $a': =a+\ell/m$ and assume 
\begin{equation}\label{encore2}
C_0 - \epsilon > a'(m + \Lambda 2^{a'-1}), 
\end{equation}
in order to use inequalities simiar to (\ref{HH})--(\ref{HH2}) with $a'$ instead
of $a$. 
Note that, at least, for $\ell>0$ -- and automatically $k$ -- small enough, the
latter inequality holds true due to 
(\ref{encore1}).
Now, let us collect all terms and their bounds, 
integrate and take expectations (exploiting 
an appropriate localization for $\tilde M$ if necessary), 
\begin{eqnarray*}
E_x L_{m,a,k}(t\wedge \tau_0,X_{t\wedge \tau_0}) 
+ (C_{m,\Lambda,a}-\epsilon) E_x \int_0^{t\wedge \tau_0} 
(1+s)^{k} L_{m,a-1/m}(X_s) \,ds 
 \\\\
\le L_{m,a}(x) 
+ C' E_x \int_{0}^{\infty} 
(1+s)^{k-1-\ell} E_x 1(s\le t\wedge\tau_0) 
L_{m,a+\ell/m}(X_s)\,ds. 
 \\\\
\le L_{m,a}(x) + C'' L_{m,a+(\ell-1)/m}(x).	\hspace{2cm}
\end{eqnarray*}
(This writing does not necessarily mean that $\ell \ge 1$.) Due to Fatou's lemma, with 
$\ell = k  + \delta$ and $\delta>0$ (i.e., $\ell > k$), this imples,   
\begin{eqnarray*}
E_x L_{m,a,k}(\tau_0,X_{\tau_0}) 
+ C' E_x \int_0^{\tau_0} 
(1+s)^{k} L_{m,a-1/m}(X_s) \,ds 
 \\\\
\le L_{m,a}(x) 
 + C'' L_{m,a+(\ell-1)/m}(x).  \hspace{2cm}
\end{eqnarray*}
Since $L_{m,a-1/m}(X_s)\ge 1$ for $s<\tau_0$ 
(notice that $a+(\ell-1)/m>0$), we get, 
\begin{eqnarray*}
\displaystyle E_x \tau_0^{k+1} 
\le C L_{m,a}(x)  + C L_{m,a+(\ell-1)/m}(x),
\end{eqnarray*}
or just 
\begin{eqnarray}\label{simply}
E_x \tau_0^{k+1} 
\le C L_{m,a+(\ell-1)_+/m}(x).
\end{eqnarray}
Notice that for $x=\Delta_0$, the inequality (\ref{simply}) also trivially
holds true. 

~

\noindent
{\bf 4.} 
The bound (\ref{simply}) along with moment inequalities (\ref{HH}--\ref{HH2})
for Markov models ``usually''
already suffice for establishing the desired rate of convergence and there are
several standard ways to accomplish the proof. So, in
principle, we may  claim our result at this point. However, we 
give a
sketch of the remaining proof 
for completeness of the presentation and to address some specifics of the
models. It is due to this specifics that while considering a couple of processes
we need to take some additional care so as to tackle the hitting time of the
``origin'' for this couple, while ``usually'' it is enough to estimate 
moments of the hitting time of some neighbourhood of the origin. 
In this {\em second part of the proof,} 
we consider two independent 
versions $X$ and $\tilde X$ of our 
Markov process, one starting at $x$ and another at the 
stationary measure $\mu$. We are going to show 
how arrange coupling. Recall that the stationary version 
exists due to the Harris-Khasminskii principle, see the remark above. 
Denote $\bar \tau_0 : = \inf(t\ge 0: X_t  = \tilde X_t = \Delta_0)$. 
Given $\tilde X_0 = y$, it may be proved that also
\begin{eqnarray}\label{taubar}
E_{x,y} \bar \tau_0^{k+1} 
\le C \bar L_{m,a+\ell/m}(x,y).
\end{eqnarray}
Indeed, let $R>0$ and for given $m$, $a$ and $\ell$, denote
\begin{eqnarray*}
\bar\tau_{0,R}: = \inf(t\ge 0: \; X_t=\Delta_0 \; 
\mbox{and}\; L_{m,a-1/m}(Y_t)\le R, 
 \\
\;\; 
\mbox{or}\;\;  Y_t=\Delta_0 \; \mbox{and}\; L_{m,a-1/m}(X_t)\le R).
\end{eqnarray*}
The idea of evaluating $\bar \tau_0$ is to establish a bound for 
$\bar\tau_{0,R}$ and then to use it for managing $\bar \tau_0$ 
with the help of (\ref{simply}). For this goal, consider Lyapunov 
functions 
$$
\bar L_{m,a}(X_{t},Y_t):= L_{m,a}(X_t) + L_{m,a}(Y_t), \;\;
\bar L_{m,a,k}(t,X_t,Y_t) 
:= (1+t)^k \bar L_{m,a}(X_{t},Y_t),
$$ 
with the same values 
of $m,a,k$ as above for a single component. We notice that 
at any moment $t$ when $(X_t, Y_t) \not = (\Delta_0, \Delta_0)$, 
the Lyapunov function $\bar L_{m,a,k}$
serves well in the sense that it decreases on average at least 
as fast as a single component one, $(1+t)^k L_{m,a}(X_t)$, say 
(if $X_t\not=\Delta_0$). 
If it occurs for the first time that $X_t=Y_t=\Delta_0$, then it means that 
$t=\bar\tau_0$. So, we have to inspect what happens at $t$ when 
$X_t=\Delta_0$ and $\Delta_0\not = Y_t$, but $L(Y_t)\ge R$, say. 
In this case the idea is that the average increment of 
$L_{m,a-1/m}(X_t)$ is, of course, positive but equals just $\lambda_0 dt$ and, hence, 
may be easily compensated by a large negative on average increment 
of the other component $L_{m,a-1/m}(Y_t)$. In this way we will establish below 
the bound
\begin{eqnarray}\label{tau0r}
 E_{x,y}\bar\tau_{0,R}^{k+1} \le C \bar L_{m,a+\ell/m}(x,y), 
\end{eqnarray}
under the condition (\ref{encore2}).
Then, once $\bar\tau_{0,R}$ occurred, we may wait some fixed time 
$t_1$ sufficient for $Y$ to achieve $\Delta_0$ with a large probability, 
say, at least $1/2$, while $X$ remains at $\Delta_0$ all that time with probability $\exp(-\lambda_0 \, t_1)$. 
If this scenario is {\em not} realised -- which probability 
does not exceed some constant $\nu<1$ 
-- then we stop at $\bar\tau_{0,R} + t_1$ or a bit earlier 
if either $X$ exits from $\Delta_0$, or $L_{m,a-1/m}(Y)$ exceeds level $R+1$ (say). 
Then we wait again until the ``next'' moment $\bar\tau_{0,R}$ 
and repeat the whole procedure of the ``attempt'' to meet 
both components at $\Delta_0$. Thus, we will evaluate $\bar\tau_{0}$ 
by means of some geometric series, which would guarantee the 
desired inequality (\ref{taubar}). Hence, let us show the bound 
(\ref{tau0r}) first. Recall that we have $C_{m,\Lambda,a}>0$ and even $C_{m,\Lambda,a'}>0$ ($a'=a+\ell/m$) due to (\ref{encore1}) and  (\ref{encore1}), and choose  
$\epsilon$ and $R$ so that   
\begin{equation}\label{eqR}
(C_{m,\Lambda,a} - \epsilon) R > \lambda_0. 
\end{equation}
Then there exists $C'>0$ such that 
$(C_{m,\Lambda,a} - \epsilon- C') R \ge  \lambda_0)$.
Denote 
\begin{eqnarray*}
e^1_t:= 1(X_t\not=\Delta_0, Y_t\not=\Delta_0), \quad
e^2_t:= 1(X_t=\Delta_0, L_{m,a-1/m}(Y_t)\ge R), 
 \\
e^3_t:= 1(Y_t=\Delta_0, L_{m,a-1/m}(X_t)\ge R). \hspace{2cm}
\end{eqnarray*}

~

\noindent
{\bf 5.} We start with the function $\bar L_{m,a}(X_t,Y_t)$ 
on $t<\bar\tau_{0,R}$. Repeating the calculus at the step 1, 
we obtain the following bounds, 
\begin{eqnarray}\label{HH3}
 E_{x,y} \bar L_{m,a} (X_{t\wedge \bar\tau_{0,R}}, 
 Y_{t\wedge \bar\tau_{0,R}}) 
+ E_x \int\limits_0^{t\wedge \bar\tau_{0,R}} 
\{[(e^1_s + e^3_s)C_{m,\Lambda,a} L_{m,a-1/m}(X_{s}) 
-e^3_s \lambda_0]
 \nonumber \\ \\ \nonumber
+ [(e^1_s + e^2_s)C_{m,\Lambda,a} L_{m,a-1/m}(Y_{s}))
- e^2_s\lambda_0]\}\, ds 
\le \bar L_{m,a} (x,y),
\hspace{1cm}
\end{eqnarray}
and, due to Fatou's lemma, also 
\begin{eqnarray}\label{HH4}
 E_{x,y} \bar L_{m,a} (X_{\bar\tau_{0,R}}, 
 Y_{\bar\tau_{0,R}}) 
+ E_x \int\limits_0^{\bar\tau_{0,R}} 
\{[(e^1_s + e^3_s)C_{m,\Lambda,a} L_{m,a-1/m}(X_{s}) 
-e^3_s \lambda_0]
 \nonumber \\ \\ \nonumber
+ [(e^1_s + e^2_s)C_{m,\Lambda,a} L_{m,a-1/m}(Y_{s}))
- e^2_s\lambda_0]\}\, ds 
\le \bar L_{m,a} (x,y).
\end{eqnarray}
Due to the condition (\ref{eqR}), all integrands ``$[\ldots]$'' in  
(\ref{HH3}) and  (\ref{HH4}) 
are non-negative, so, in particular, for any $t\ge 0$,
\begin{eqnarray}\label{eleqe}
E_{x,y} \bar L_{m,a} (X_{t\wedge \bar\tau_{0,R}}, 
 Y_{t\wedge \bar\tau_{0,R}}) \vee 
E_{x,y} \bar L_{m,a} (X_{\bar\tau_{0,R}}, 
 Y_{\bar\tau_{0,R}}) \le \bar L_{m,a} (x,y).
\end{eqnarray}

~

\noindent
{\bf 6.} Now we are ready to consider the Lyapunov function $\bar L_{m,a,k}(t,X_t,Y_t)$ depending also on time. 
Similarly to the step 3 -- see the formula (\ref{twoterms}) --  
we have on $t<\bar \tau_{0,R}$ with some new local martingale $\hat M_t$, 
\begin{eqnarray*}
d\bar L_{m,a,k}(t,X_t,Y_t) 
= \bar L_{m,a,k}(t+dt, X_{t+dt}, Y_{t+dt}) - 
\bar L_{m,a,k}(t,X_{t},Y_t) 
 \\\\
\le 
e^1_t (1+t)^{k} 
\left[(I^Y_1 - I^Y_2 + I^Y_3 + \frac{k L_{m,a}(Y_t)}{1+t} )\,dt \hspace{2cm}
 \right.\\\\\left.
+ (I^Y_1 - I^Y_2 + I^Y_3 + \frac{k L_{m,a}(X_t)}{1+t})\,dt 
+ d\hat M_t\right]  \hspace{2,5cm} 
 \\\\
+ e^2_t  (1+t)^{k} 
\left[(I^Y_1 - I^Y_2 + I^Y_3 
+ \lambda_0 + \frac{k L_{m,a}(Y_t)}{1+t}) \,dt 
+ d\hat M_t\right]  \hspace{1cm}
 \\\\
+ e^3_t  (1+t)^{k} 
\left[(I^X_1 - I^X_2 + I^X_3 
+ \lambda_0 + \frac{k L_{m,a}(X_t)}{1+t})\,dt +  d\hat M_t\right] 
 \\\\
=: (J_1 + J_2 + J_3)dt + d\tilde M_t,   \hspace{3cm}
\end{eqnarray*}
again with a new local martingale $\tilde M$ and with
\[
J_1 := 
e^1_t (1+t)^{k} 
\left[I^Y_1 - I^Y_2 + I^Y_3 + \frac{k L_{m,a}(Y_t)}{1+t}  
+ I^Y_1 - I^Y_2 + I^Y_3 + \frac{k L_{m,a}(X_t)}{1+t}\right], 
\]
\[
J_2 := 
e^2_t  (1+t)^{k} 
\left[I^Y_1 - I^Y_2 + I^Y_3 
+ \lambda_0 + \frac{k L_{m,a}(Y_t)}{1+t}
\right],
\]
\[
J_3 := 
e^3_t  (1+t)^{k} 
\left[I^X_1 - I^X_2 + I^X_3 
+ \lambda_0 + \frac{k L_{m,a}(X_t)}{1+t}\right].
\]
Here the first term $J_1$ is estimated identically to what was 
done at the step 3  for the only component $X$, and this gives 
\begin{eqnarray*}
J_1 
\le  e^1_t\left(-(1+t)^k (C_{m,\Lambda,a} - \epsilon) 
\bar L_{m,a-1/m}(X_{t},Y_t)\right)
 \\ 
+ e^1_t k(1+t)^{k-1} \,  
\times \frac{k^\ell}{(\epsilon (1+t))^\ell}
\bar L_{m,a+\ell/m}(X_t, Y_t)^{}.
\end{eqnarray*}
The second and the third terms allow the bounds, 
\begin{eqnarray*}
J_2
\le  e^2_t\left(-(1+t)^k C_{m,\Lambda,a} 
L_{m,a-1/m}(Y_{t}) - \lambda_0 \right)
 \\
+ e^2_t k(1+t)^{k-1} \,  
\times \frac{k^\ell}{(\epsilon (1+t))^\ell}
L_{m,a+\ell/m}(Y_t)^{},
\end{eqnarray*}
\begin{eqnarray*}
J_3
\le  e^3_t\left(-(1+t)^k C_{m,\Lambda,a} 
L_{m,a-1/m}(X_{t}) - \lambda_0 \right)
 \\ 
+ e^3_t k(1+t)^{k-1} \,  
\times \frac{k^\ell}{(\epsilon (1+t))^\ell}
L_{m,a+\ell/m}(X_t)^{},
\end{eqnarray*}

Now, let us collect all terms and their bounds, 
integrate and take expectation, also using 
localization for the martingale term if necessary. Notice that 
$1(s<\bar\tau_0)(e^1_s+e^2_s+e^3_s)=1(s<\bar\tau_0)$ 
and 
\begin{eqnarray*}
1(s<\bar\tau_0)[(e^1_s + e^3_s) L_{m,a-1/m}(X_s) 
+ (e^1_s + e^2_s) L_{m,a-1/m}(Y_s)] 
 \\
= 1(s<\bar\tau_0) \bar L_{m,a-1/m}(X_s,Y_s). \hspace{2cm}
\end{eqnarray*}
So, we have,  
\begin{eqnarray*}
E_{x,y} \bar L_{m,a,k}(t\wedge \bar\tau_{0,R}, 
X_{t\wedge \tau_0}, Y_{t\wedge \tau_0}) \hspace{2.5cm}
 \\\\
+ (C_{m,\Lambda,a}-\epsilon-\frac{\lambda_0}{R}) 
E_{x,y} \int_0^{t\wedge \bar\tau_{0,R}} 
(1+s)^{k} \bar L_{m,a-1/m}(X_s,Y_s) \,ds 
 \\\\
\le \bar L_{m,a}(x,y) 
+ C' \int_{0}^{\infty} 
E_{x,y} 1(s\le t\wedge\bar\tau_{0,R}) 
(1+s)^{k-1-\ell} \bar L_{m,a+\ell/m}(X_s,Y_s)\,ds. 
\end{eqnarray*}
Further, recall that $a'=a+\ell/m$ and (\ref{encore2}) holds true, whence, 
\begin{eqnarray*}
E_{x,y} 1(s\le t\wedge\bar\tau_{0,R}) 
\bar L_{m,a+\ell/m}(X_s,Y_s)
\le \bar L_{m,a+\ell/m}(x,y).
\end{eqnarray*}
From here we conclude,    
\begin{eqnarray*}
E_{x,y} \bar L_{m,a,k}(t\wedge \bar\tau_{0,R}, 
X_{t\wedge \tau_0}, Y_{t\wedge \tau_0}) 
+ C' E_{x,y} \int_0^{t\wedge\bar\tau_{0,R}} 
(1+s)^{k} \bar L_{m,a-1/m}(X_s,Y_s) \,ds 
\\\\
\le C \bar L_{m,a+\ell/m}(x,y).\hspace{3.5cm}
\end{eqnarray*}
Due to Fatous's lemma, with 
$\ell = k  + \delta$ (i.e., $\ell > k$) this imples,   
\begin{eqnarray*}
E_{x,y} \bar L_{m,a,k}(\bar\tau_{0,R},X_{\bar\tau_{0,R}},Y_{\bar\tau_{0,R}}) 
+ C' E_{x,y} \int_0^{\bar\tau_{0,R}} 
(1+s)^{k} \bar L_{m,a-1/m}(X_s,Y_s) \,ds 
 \\\\
\le C \bar L_{m,a+\ell/m}(x,y). \hspace{2.5cm}
\end{eqnarray*}
Since $\bar L_{m,a-1/m}(X_s)\ge 1$ on $s<\bar \tau_0$ (and on $s<\bar \tau_{0,R}$), 
we get 
\begin{eqnarray*}
\displaystyle E_{x,y} \bar\tau_{0,R}^{k+1} 
\le C \bar L_{m,a+\ell/m}(x,y),
\end{eqnarray*}
so that (\ref{tau0r}) is established. 

~

\noindent
{\bf 7.} 
Now let us show (\ref{taubar}).
As explained above, to this aim we choose $t_1$ so that
$$
\sup_{u: \, L_{m,a-1/m}(u)\le R+1} \, 
P_{u}(\tau_0>t_1) \le 
t_1^{-(k+1)} \sup_{u: \, L_{m,a-1/m}(u)\le R+1}\, E_{u}\tau_0^{k+1} \le \frac12. 
$$
Recall that $\nu<1$ is defined above in the step 4 as follows, 
\begin{eqnarray*}
1-\nu:= \inf_{L_{m,a-1/m}(u)\le R+1} P_{\Delta_0, u}\left(X_s\equiv \Delta_0, \, 
0\le s\le t_1, \, \& \, 
\exists \, t \in [0, t_1]: \, Y_t = \Delta_0\right)
 \\\\
\ge \frac12\, \exp(-\lambda_0 \, t_1) > 0. \hspace{3.5cm}
\end{eqnarray*}
Let 
\(
\bar\tau(1):= \bar\tau_{0,R}, \;\; 
\bar\tau(n+1):= \theta_{\bar\tau(n)+t_1}\bar\tau_{0,R}, 
\;\; n=1,2,\ldots, 
\)
where $\theta_t$ is a shift operator for the process 
$((X_t,Y_t), \, t\ge 0)$ 
(see \cite{Dynkin}). 
Then we estimate, 
\begin{eqnarray}\label{sumtau0}
 E_{x,y}\bar \tau_{0}^{k+1} \le \sum_{n\ge 1}^{} 
(E_{x,y} (\bar\tau(n) + t_1)^{k+1} \nu^{n-1}(1-\nu).
\end{eqnarray}
Whence,   
\begin{eqnarray*}
E_{x,y} (\bar\tau(1) + t_1)^{k+1} 
\le 2^k \,(E_{x,y} \bar\tau(1)^{k+1} + t_1^{k+1})
\le C \bar L_{m,a+\ell/m}(x,y) + C. 
\end{eqnarray*}
By induction and using 
$\bar\tau(n) 
= \sum_{k=1}^{n}(\bar\tau(k) - \bar\tau(k-1))$ 
with $\bar\tau(0):=0$, we have, 
\begin{eqnarray*}
E_{x,y} (\bar\tau(n) + t_1)^{k+1}
\le C n^{k}((n-1) +  
\bar L_{m,a+\ell/m}(x,y) + 1)
 \\
\le Cn^{k+1}(\bar 
L_{m,a+\ell/m}(x,y) + 1). \hspace{2cm}
\end{eqnarray*}
Substitution of this into (\ref{sumtau0}) shows that, 
indeed, (\ref{tau0r}) implies (\ref{taubar}), under the assumption (\ref{eqR}), the latter being guaranteed by (\ref{encore2}).

~

\noindent
{\bf 8.} 
From (\ref{taubar}) we conclude, 
with invariant distribution $\mu$,
\begin{eqnarray}\label{taubar3}
E_{x,\mu} \bar \tau_0^{k+1} 
\le C L_{m,a+\ell/m}(x) 
+ C \int L_{m,a+\ell/m}(y)\,\mu(dy).
\end{eqnarray}
Recall that from (\ref{HH})  it follows that $\mu$ integrates the function 
$L_{m,a-1/m}(x)$ for any couple $(a,m)$ satisfying 
$a, m > 1$ and (\ref{am}): $C_0> a(m + \Lambda 2^{a-1})$.  
Hence, for the integral in (\ref{taubar3}) to converge, it suffices  
\(
C_0> (a + \ell/m)(m + \Lambda 2^{a -1 + \ell/m}). 
\)
The latter is guaranteed by  (\ref{encore2}), that is, by (\ref{encore1}) and by the choice of $\ell$ close enough to $k$.
In turn, (\ref{encore1}) is guaranteed by the assumption 
(\ref{aslambda}) if $k>0$ and $\ell>0$ are sufficiently small. Then, for any $k>0$ 
small enough, there exist $a>1, \, m>1$ and $\ell > k$ such that 
the integral in (\ref{taubar3}) converges and we get 
\begin{eqnarray}\label{taubar4}
E_{x,\mu} \bar \tau_0^{k+1} 
\le C L_{m,a+\ell/m}(x) 
+ C.
\end{eqnarray}

~

\noindent
{\bf 9.} 
Now, we may estimate 
the right hand side in the coupling inequality,  
$$
\|\mu^x_t - \mu\|_{TV} = 2\sup_A(\mu^x_t - \mu)(A) \le 
2 P_{x,\mu}(T>t), 
$$ 
where $T:= \inf(t\ge0: \, X_t=\tilde X_t = 0)$. It follows from 
(\ref{taubar4}) in a standard way (cf. \cite{Kalashnikov93}, \cite{Ve00}, et al.) 
that for any $\nu>0$ there exists $C>0$ such that 
\begin{equation}\label{lastineq}
 P_{x,\mu}(T>t) \le 
 C (1+L_{m,a+\ell/m}(x))(1+t)^{-(k+1)+\nu}.
\end{equation}
This is equivalent to (\ref{main}). 
The Theorem \ref{Thm1} is proved.

~

\noindent 
{\it Remark 5.} The main result may be extended to 
the case where both $\lambda$ 
and $h$ depend 
on the whole state of the process $X$ satisfying 
the same generic assumptions (\ref{ash}), (\ref{aslambda}) 
and (\ref{Lambda}), with $h(t)$ replaced by 
$h(t,x)$ and $\Lambda:=\sup_{n\ge 1}(\lambda_n/n)$ by 
$\Lambda:=\sup_{n\ge 1,\,x}(\lambda_{n, x}/n)$. 
Similar convergence rate {\em independent of $N$} may be proved
in the same way for the model with any $N<\infty$; here ``usual'' bounds could
easily depend on this parameter. 
Similar bounds may be established for 
{\em mixing rates} 
by using the approach from \cite{Ve00}. 
For a {\em random} initial distribution $\mu_0$, 
similar or weaker bounds may be proved depending on moments of $\mu_0$. 

\section*{Acknowledgements}
The author is grateful to G. A. Zverkina and to two anonymous Referees for many useful remarks.

%




\end{document}